\documentclass{article}
\usepackage{amsmath}
\input amssymb.sty
\newtheorem{prop}{Proposition}[section]

\newtheorem{theo}[prop]{Theorem}

\let\we=\wedge

\def\ra{{\rightarrow}}
\def\a{{\alpha}}
\def\b{{\beta}}
\def\e{{\epsilon}}
\def\nn{{\noindent}}

\def \N{{\mathbb  N}}

\def \R{{\mathbb R}}
\def\l{{\lambda}}
\def\ts{\times}
\let\l=\lambda

\title{Timescales of population rarity and commonness in 
random environments  }
\author{REGIS FERRIERE\thanks{Unit of Mathematical Evolutionary Biology, Department of Ecology, CNRS UMR
7625, \'Ecole Normale Sup\'erieure, 46 rue d'Ulm, 75230 Paris cedex 05, France
(e-mail: ferriere@biologie.ens.fr)
2 Department of Ecology and Evolutionary Biology, University of Arizona, Tucson AZ
85721, USA};
\,\,
ALICE GUIONNET\thanks{\'Ecole Normale Sup\'erieure de Lyon,
Unite de Math\'ematiques pures et appliqu\'ees,
UMR 5669,
46 All\'ee d'Italie,
69364 Lyon Cedex 07, France},\
\,\,
IRINA KURKOVA
\thanks{Laboratoire de Probabilit\'es et Mod\`eles Al\'eatoires,
  Universit\'e Paris 6, B.C. 188; 4, place Jussieu, 75252 Paris Cedex
  05, France. }}
\date{22/09/2004}
\def\d{{\delta}}
\begin{document}

\maketitle

\nn {\bf{Abstract:}}
This paper investigates the
influence of environmental noise on the characteristic
 timescale of the dynamics of density-dependent populations. 
General results 
are obtained on  the statistics of time spent in rarity (i.e.\
 below 
a small threshold on
population density) and time spent in commonness (i.e. above a
 large threshold). The nonlinear stochastic
 models under consideration form
 a class of Markov chains on the state space $]0,
\infty[$ which are transient 
if the intrinsic growth
rate is negative and  recurrent if it is positive
or null. 
 In the
recurrent  case,  we obtain a  necessary and sufficient 
   condition for positive recurrence and precise estimates 
 for the distribution of times
of rarity and commonness. In the null recurrent, critical
 case  that applies to ecologically
neutral species, the distribution of rarity time is a 
universal power law with
exponent $-3/2$. This has
implications for our understanding of the long-term dynamics of some natural
populations, and provides a rigorous basis for the
 statistical description of on-off
intermittency known in physical sciences.

Key words or phrases:  Population dynamics - 
  Stochastic non-linear
difference equations 
  - Environmental
stochasticity - Rarity - Ecological timescales - Power law -
 On-off intermittency - Markov chains - Martingales.

\section{Introduction}

Many of the traditional approaches in ecological theory are based
 on a paradigm that the
ecological systems that we observe in nature correspond in some 
way to stable states of
relatively simple ecological models. This is a view embodied in 
many classic approaches
to ecological problems, ranging from the use of models to look at
 the coexistence of
species to the interpretation of short-term experiments [May
\cite{M1}
]. This approach,
however, has been questioned. 
As ecologists emphasize that long-term experiments yield
different results compared with short-term experiments 
[Brown et al. \cite{R}], theory must
take the role of timescales into account and look beyond an emphasis on determining
asymptotic stable states [Hastings \cite{H}].

 Of particular applied interest is the issue of
identifying timescales that are relevant for population management. Study of pests or
disease outbreaks often have the time spent by the system in quiescence as one of their
focus. The conservation of threatened species requires an understanding of how
individual and environmental parameters affect the time spent by a population above
some predefined critical threshold of abundance. Planning 
the exploitation of renewable
resources in an uncertain world may also benefit from
 probabilistic estimates of time spent below
and above predefined thresholds of ecological
and economic  interest.
 This points to a common mathematical problem: given
information about the current population state, how can one 
characterize
 the distribution of time
spent (or needed) by a population trajectory 
within (or to reach) a particular domain of
the population state space?
 Here we address this problem for a class of unstructured,
density-dependent stochastic population
 models which includes the well-known
Ricker model.

Theoretical ecologists have long used simple, unstructured models
to investigate general properties of population dynamics. These
models, however, often miss essential features of natural
populations, yielding predictions that turn out to be wrong even
qualitatively [Durrett and Levin 1994 \cite{DL}]. The difference equation
introduced by Ricker \cite{R}  makes a remarkable
exception. Although it was originally introduced in the tradition
of phenomenological modelling, the Ricker model has recently been
re-derived from first principles accounting for the discreteness
of individuals, the stochastic nature of birth and death events,
and the spatial localization of interactions between individuals
[Royama \cite{Ro} , Ripa and Lundberg \cite{RL}, 
Sumpter and Broomhead \cite{SB}].
In a constant environment, the model reads
\begin{equation}\label{eq1}
 X_{n+1} = X_n
f(X_n)\end{equation} where 
\begin{equation} 
\label{fff}
f(X_n) = \exp[r-aX_n ],
\end{equation}
 and $X_0 >
0$. $X_n$ measures the population size (or density at time $n$).
The number $r$ measures the per capita growth rate at low
population density when density-dependent effects are negligible;
thus, $\l=\exp(r)$
 is the population's
geometric rate of increase from rarity.  Unlimited growth is prevented
by the term $\exp(-a X_n)$, where ${\it a}>0$ measures the intensity of
negative interactions between individuals. This non-linear term
typifies the mode of density dependence called overcompensatory,
whereby ${\it xf(x)}$ decreases toward zero as ${\it x}
$ becomes larger than some
threshold (with undercompensatory density dependence, $xf(x)$ 
would be
monotonically increasing). 

\smallskip 

\noindent{\bf Notation.} Hereafter
we set $\l=f(0)$ for the model (\ref{eq1})
  with general fucntion $f(x)>0$; $\l$ is called the intrinsic
growth rate.

\smallskip

The Ricker model shows a somewhat intricate transition from
equilibrium and cycles to chaos as $r$ increases [May \cite{M2}]. 
The natural timescale of the Ricker model in a constant
environment is essentially determined by the intrinsic growth
rate. When $\l<1$, the origin is the only stable equilibrium state
of the system (hence $\ln X_n$ goes to $-\infty$ and is thus
transient), whereas when $\l>1$, there are two equilibria : the
origin, always unstable; and a positive equilibrium which can be
stable or not depending on the specific value of $\l$ (hence $\ln
X_n$ is recurrent). Multi-dimensional 
Ricker models are also of interest to
describe interacting populations
of different species \cite{hhj}.

 Randomness, however, is inherent to the real
world. Perhaps the simplest way 
of incorporating  stochasticity  into
(\ref{eq1}) (cf. \cite{Ni}) is to set
\begin{equation} 
\label{eed}
 X_{n+1} = X_n f(X_n) \exp(Y_n)
\end{equation}
where $Y_n, n=0, 1, 2,\cdots$
 are i.i.d random variables with zero expectation; $(X_n)$
 is then a
Markov chain (MC) on the state space $]0, \infty[$. This noisy
Ricker model was first introduced in [Kornadt et al. 1991] and
considered in the context of MC theory by [Gyllenberg
and all \cite{GM1,GM2}]. The random variable $Y_n$
 can be seen as an additive
perturbation of the intrinsic growth rate, thus providing a
model of environmental stochasticity \cite{GM2};  or as a
perturbation of the population density which is `felt'
effectively  by each individual in the population--a form of
stochasticity that we call {\it random heterogeneity}. Random 
heterogeneity arises from differences in competitive 
abilities among individuals which may be rooted in
intrinsic individual differences or in spatial variation in
habitat quality \cite{ho}. 
In the first case a simple choice is to assume the
$Y_n$ to be Gaussian whereas in  the latter case, a typical choice for
the law of $Y_n$ is the log-normal distribution \cite{Ni}.
This model  generalizes 
straightforwardly  to the case of $N$ species in interaction.

\medskip 

\def\e{{\epsilon}}

\medskip

  Our goal in this paper is to study the time spent by the
population within given ranges of density. Of particular interest
are the time $T_\epsilon $ spent in a state of rarity defined by a
small density threshold $\e>0$,
\begin{equation}\label{eq13}
T_\epsilon:=\inf\{n\ge 1 : X_n\ge \e\}
\end{equation}
when the current density $X_0$ is smaller than
$\epsilon$,   and the
time $T_M$ spent in a state of commonness defined by a large
density threshold $M>0$, 
\begin{equation}\label{eq14}
T_M:=\inf\{n\ge 1 : X_n\le M\}
\end{equation}
 when the current density $X_0$ is larger than $M$.
 More generally, we define
the exit time of the Markov chain $(X_n)$
 starting at some $ X_0
> 0$ from a domain $A$ of population densities (in $\R^+$) as
\begin{equation}\label{eq15}
T_A=\inf\{n\ge 1 ; X_n\in A^c\}
\end{equation}
 so that $T_\epsilon = T_{]0, \epsilon[}$ and $T_M
= T_{]M, \infty[}$. We further define the time
of medium abundance as $T_{[\e,M]}$,
and the time of escape from extremes as $T_{[\e,M]^c}$.
We then  ask the following
questions. 

\begin{itemize} 

\item What is the qualitative long-term 
pattern of population dynamics ?
In mathematical terms, we ask whether 
the MC is recurrent, null recurrent
or  transient. This question was already adressed in \cite{FH}
under the assumptions 
 that the law of the random noise
has exponential moments. Here we relax this assumption 
to consider moments of order strictly greater than one, thereby
 extending
the scope of the model to random heterogeneity.

\item  Then, we
seek general results on the tail
of the distribution of times of rarity
or commonness which 
are independent of the current population
density $X_0$. This means asking 
how the probabilities  $P(T_\epsilon > n)$ (if $X_0<\epsilon$), 
 $P(T_M > n)$ (if $X_0>M$)
and   $P(T_{[\epsilon,M]^c} > n)$ (if $X_0\not\in [\e, M]$) 
 behave when ${\it n}$ goes to infinity. Almost equivalent 
is to ask for which 
 $p>0$ is  $E[(T_A)^p]$ finite, given $A=]0,\epsilon[$,  
$A=]M,\infty[$ or $A=]0, \epsilon[\cup ]M,\infty[$
with $X_0 \not\in A$.

\item  Another interesting  feature of  the 
distribution of the times of rarity 
or commonness concerns the dependence of
their moments on the current density $X_0$.
Therefore, we shall study 
 the behavior  
   of the $p$th moment
    of $T_A$, i.e.\ $E_{X_0} (T_A)^p$ for given $A$ and $p$, 
    as a function of the initial state 
  $X_0$; namely,  for given $A=]0,\e[$   when 
   $X_0\to 0$, or  for given $A=]M, \infty[$  when  $X_0\to \infty$.

\item How much do the preceeding results depend on the
choice of the function~$f$ and the distribution of the noise ? 

\item How do the results extend to multispecies interactions ? 
\end{itemize}

We shall consider the noisy 
dynamical system (\ref{eed}) with a general continuous
function $f:\R^+\ra\R^+$  
and independent copies  $(Y_i,i\ge 0)$
of a variable $Y$ such that $E[Y]=0$.
Throughout this article we
shall assume the following :

\begin{enumerate}
\item $E[|Y|^{1+\d}]<\infty$ for some $\d>0$.
\item 
The function $f(x)>0$ is continuous 
on $[0, \infty[$ and there exist $a>0$, $r\in \R$ 
    such that
  \begin{equation}
\label{ara}
\lim_{x\to \infty} f(x)e^{ax-r}=1.
\end{equation} 
\end{enumerate}

The latter assumption says that at
large density population growth is well
approximated by the Ricker model.
A crucial consequence is that the results 
will then  only depend
on whether the parameter $\l$
is strictly greater, smaller or equal to one,
that is, whether the population
is intrinsically growing, or declining, or neutral. 
   For the Ricker model specified by (\ref{fff}), 
  these conditions translate into $r<0$, $r>0$ and $r=0$.

We recall the meaning of basic 
 terminology used for the
classification  MC processes. 
   We say that the MC is {\it transient\/}
  if there is a positive probability that
  the time taken by the  process, when initiated in any bounded segment 
  $A=[a,b]\subset ]0, \infty[$ (with $0<a<b<\infty$),
 to return to $A$, is
infinite; in other words, the probability that the MC never returns to
 $A$ is positive. The
MC is {\it recurrent\/} if the return time 
to $A$ is finite almost surely.
 The MC is {\it positive
recurrent} if, in addition, the expectation of the return time to $A$
 is finite too 
 and {\it null recurrent\/} if the expectation of this
 time is infinite. 

 To help developing biological intuition for
our analysis, it may be worth emphasizing that,
heuristically, 
 a (single) population's transience 
will indicate either extinction
or escape toward infinite density. Finding that
the distribution of time spent in a certain state 
has an exponential tail indicates that
the population typically visits that
states for short time only; in contrast,
a `polynomial' or `heavy' tail tells that
the population resides very long in that state.

 Recurrent notations in this paper
will be  $P_x$ for the law of $(X_n, n\ge 0)$
 starting from $X_0 = x$, and $E_x$ for the
expectation under $P_x$.

\section{Time of commonness and time of medium abundance}

The distribution of the random time  $T_M$ of commonness 
(if $M$ is fixed large enough and $X_0>M$)  
    is described by  the following  Theorem~\ref{th0}. 

\begin{theo}\label{th0}

\noindent (a) {\it {\bf Asymptotics when $X_0 \to \infty$.} }
     For  any   $X_0>M>0$
  we have:  $E_{X_0} T_M<\infty$.
  Furthermore
\begin{equation}
\label{f1}
 E_{X_0}T_M \to 1, \ \  \hbox{ as }X_0 \to \infty,
\end{equation}
   and  $T_M\to  1$ a.s.
  as $X_0 \to \infty$.

\noindent (b) {\it\bf{ Asymptotics of $P_{X_0}(T_M>n)$
   as $n\to \infty$, when  $X_0$ is  fixed. }}

  From  the fact that $E_{X_0} T_M< \infty$, we have
   $P_{X_0}(T_M>n)=o(n^{-1}\ln^{-1}n)$
 as $n\to \infty$.

 Assume that in addition  for some  $\a_0>0$
$E[e^{\a_0 Y}]<\infty$. Then
there exists $\kappa<1$
such that for any $X_0>M$ and any $n\geq 1$
\begin{equation}
\label{e1}
P_{X_0} (T_M\ge n)\le C(X_0)\kappa^n
\end{equation}
    with some constant $C(X_0)$.
 Moreover for all $M$ large enough
     $$\kappa\leq \inf_{\a>0} \sup_{x\ge M} f(x)^\a
E[ e^{\a Y}]<1.$$ 
 If this infimum is reached for $\a=\a_0$, then
 $C(X_0)\leq (X_0/M)^{\a_0}$.

\end{theo}
Part (a) of the theorem says that the
negative density-dependence can pull the population 
out of extreme commonness in no longer than
one step. Moreover, the distribution
of the time of commonness
    is qualitatively little sensitive to the intrinsic
growth $\l$:
   its expectation $E_{X_0}T_M$ 
      is always finite.
 Consequently, the tail of the distribution
  of $T_M$,
 $P_{X_0}(T_M>n)$ decreases faster
   than $(n\ln n)^{-1}$ as $n\to \infty$ for any
current density $X_0$
larger than $M$.
 In the case where 
the random noise has exponential moments
($E[e^{\a Y}]<\infty$),  this
tail
even decreases exponentially, as shown by equation
(\ref{e1}))
 where the  constant $0<\kappa<1$ is independent from
the current density  $X_0$, 
    and the function  $C(X_0)$ 
 can be computed from
the law of the environmental noise $Y$.
Such estimate comes from the fact that 
the event of a large time of commonness commonness 
occurs only when an event of probability strictly less than one
is repeated  a number of times proportional to $n$ (which 
is a standard large deviation estimate).

The distribution of the random time
of medium abundance   $T_{[\e,M]}$ has 
also an exponential tail.

\begin{theo}\label{th02}
Assume that $Y$ is not uniformly bounded almost surely,
i.e for all $-\infty<-M'<0<M<\infty$, $P(-M'\le Y\le M)<1$.
Then, for all $M\ge \e>0$, there exists
$\kappa(\e,M)<1$ such that for
all $X_0\in [\e,M]$, 
$$P_{X_0}(T_{[\e,M]}\ge n)\le \kappa(\e,M)^n.$$
\end{theo}
The hypothesis that the environmental noise
 $Y$ is not uniformly bounded is 
satisfied by common choices 
of laws such as the normal or
the log-normal laws, so that we did not try
to weaken it.

Let us now consider the two other random times
$T_\epsilon$
and $T_{[\e,M]^c}$, the  behavior of which 
depends on whether  $ \l < 1, \l > 1,$ and $\l = 1$.

\section{Intrinsically declining populations}

When $\l < 1$, the MC
 is transient with  escape to $0$.
    Thus, for any rarity threshold $\epsilon>0$ 
and
 any current density smaller than $\e$,  there is a positive
probability that the time $T_\e$ 
 needed to cross the rarity threshold $\epsilon$
is infinite (\ref{t1}).  The
behavior of the process is then trivial, as it follows from
the law of large numbers \cite{FH}. The following formal
statement is given only for the sake of completness.

\begin{theo}
\label{th3}
 Assume $\l<1$ and $E[|Y|]<\infty$.

Then the MC is transient. This means that for any compact
 subset $A=[a, b]\subset ]0,\infty[$ and any $X_0\not\in A$
  $P_{X_0}(T_{A^c}=\infty)>0$.

Furthermore for any $\epsilon>0$ and any $X_0<\epsilon$,
\begin{equation}
\label{t1}
 P_{X_0}(T_\e=\infty)>0.
\end{equation}
\end{theo}

\section{Intrinsically  growing populations}
When $\l > 1$, the population MC
 is recurrent. Consider first the timescale of rarity:
a small rarity threshold 
  $\epsilon>0$ is fixed and the current density 
  $X_0$ is smaller than $\e$.
The expected time of rarity, $
E_{X_0} (T_\epsilon)$, is always finite. 
  Furthermore, we show that for any $X_0$ small enough,
 the expected time of rarity remains
of the order of minus the logarithm of initial density, i.e.
 there are two positive constants
$c_1=c_1(\epsilon)>0$ and $c_2=c_2(\epsilon)>0$ such that
$c_1 < E_{X_0} (T_\epsilon) / |\ln X_0| < c_2$
   for any initial current density 
  $X_0$ smaller than $\e$. 
Similar to the case of the time of commonness (c.f Theorem \ref{th0}),
the statistics of rarity time are given by
$ P_{X_0}(T_\epsilon > n) $ going to zero
faster than  $ (n
(\ln n))^{-1}$ when considering very
long time $n$,
and this result holds irrespective of the current density
$X_0$.
 Under some additional conditions on the environmental noise,
namely $E[e^{(-\a Y)}] < \infty$ for some $\a>0$,
 the distribution  of time spent in rarity
   has an exponentially decaying tail. In other words 
it  is controlled
according to $P_{X_0}(T_\epsilon > n) < C(X_0) \rho^n$ for any
time
$n \ge 1$; here $\rho < 1$
 is a constant
independent of the current density $X_0$, whereas $C(X_0)$
 depends upon  $X_0$.

\begin{theo}
\label{th11} Assume $\l>1$.
Then, the  MC is recurrent, i.e for any
bounded segment $A=[a,b]\subset ]0,\infty[$, $T_{A^c}$
is almost surely finite: 
  $P_{X_0}(T_{A^c}<\infty)=1$ for any $X_0>0$. 
 It is positive recurrent, 
 i.e. $E_{X_0} T_{A^c}<\infty$ for any $X_0>0$, iff 
$E [e^{Y}]<\infty$.
 For any $X_0<\e$ we have $E_{X_0} T_{\e}<\infty$, so
$P_{X_0}(T_\e>n)=o((n\ln n)^{-1})$ as $n\to \infty$. We  also have
\begin{equation}
\label{f2}
E_{X_0}T_{\e}=O(|\ln X_0|), \ \ X_0 \to 0.
\end{equation}
Moreover this order is precise in the sense that :
\begin{equation}
\label{f3}
\lim \inf_{X_0 \to 0} \frac{E_{X_0} T_{\e}}{|\ln X_0|}> 0.
\end{equation}
Assume that for some $\a>0$ $Ee^{-\a Y}<\infty$.
Then there exists $\rho<1$
such that for any $0<X_0<\epsilon$ and any $n\geq 1$
\begin{equation}
\label{e2}
P_{X_0} (T_\e\ge n)\le C(X_0)\rho^n
\end{equation}
    with some constant $C(X_0)$.
 Moreover, $\rho \leq \inf_{\a>0}e^{-\a \inf_{[0,\e]}\ln f(x)} Ee^{-\a Y}<1$.
 If this infimum is reached for $\a=\a_0$, then
 $C(X_0)\leq (\e/X_0)^{\a_0}$.

\end{theo}

Finally, we obtain
information  about the time $T_{[\e, M]^c}$ 
of escape from extremes, that is the time
  needed for
the population to reach  a
range of density $[\e, M]$ when initiated either below 
the rarity threshold 
$\e$,
 or above the commonness  threshold $M$. In
particular, we compute precise estimates for the latter's law, from
which
 we derive the
asymptotics of $P_{X_0}(T_{[\epsilon, M]^c} > n)$ 
for long times $n$. 
We show that this tail is 
of polynomial order $\sim n^{-\a}$  with degree
given by the largest  number $\a>0$ such that
$E[e^{\a Y}]$ is finite; thus, the tail depends on
the law of $Y$ via the value of $\a$. This heavy tail arises 
 from
the very long time spent by the population regrowing from
low density every time it falls short of passing below
the rarity threshold  $\epsilon$.

\begin{theo}
\label{th12} Assume $\l>1$.
For any $aM>r>\e>0$,  any $X_0 \not\in [\e, M]$
          and any $\a \in \R^+,$ 
\begin{equation}
\label{f6}
E[e^{\a Y}]<\infty\quad \Rightarrow \quad
E_{X_0}[T_{[\e, M]^c}^\a (\ln T_{[\e,
M]^c})^{(\a-1-\eta)\wedge 0}]<\infty\ 
\ \forall \eta>0
\end{equation}
   and conversely 
\begin{equation}
\label{f5}
E[e^{\a Y}]=+\infty \quad \Rightarrow \quad
E_{X_0}[T_{[\e, M]^c}^{\a}] =+\infty.
\end{equation}
In particular, if we let
 $\a_0=\sup\{\a\in \R^+ :
E[e^{\a Y}]<\infty\}$, then for any $\eta>0$
$$\liminf_{n\ra\infty} n^{\a_0-\eta+1} P_{X_0}(T_{[\e,M]^c}=n)=0$$
but
$$\limsup_{n\ra\infty} n^{\a_0+\eta+1}P_{X_0}(T_{[\e,M]^c}=n)=+\infty.$$

\end{theo}

This result calls for
two further comments. 
First, compared to \cite{FH}, we show that the chain is
positive recurrent if $\l>1$ 
  under the sole condition that $E[|Y|^{1+\d}]<\infty$
(whereas \cite{FH} assumed that
 $E [e^{\a Y}]<\infty$ for some $\a>0$) and also specify
the tail of $T_{[\e,M]^c}$.
Second,  in (\ref{f3}), we in fact show (see (\ref{gozo})
and (\ref{gozo2})) that there exists $c,c'>0$ such that
for sufficiently small $X_0$
\begin{equation}\label{eqc}
c'(\sup_{x\le\e}\ln f(x))^{-1}\le 
\frac{ E_{X_0}[T_\e]}{\ln\e-\ln X_0}\le c(\ln f(0))^{-1}=c(\ln
\l)^{-1}.
\end{equation}
This already indicates 
that the expectation of $T_\e$
is a function of $\l$ and $\e$ 
which  goes to infinity when 
 $ \l$ goes to one
and $\e$ goes to zero,
as stated in Theorem~\ref{th4}. 

\section{Intrinsically neutral populations }
 When $\l = 1$, we assume that $E[|Y|^{2+\delta_1}]<\infty$ 
  for some $\delta_1>0$.
We also need to assume that $f(x) \to 1$ as $x\to 0$ not too
  slowly: more precisely $f(x)=1+o(|\ln x|^{-1-\delta_2})$ 
  as $x\to 0$ with some small $\delta_2>0$. 
This means that the population is neutral or
almost neutral over a `fair' range of small 
densities.
If $f$ were to go to one extremely slowly, one can intuitively 
guess that the effect of neutrality 
would  be offset by the tendency of the population
to grow or decline at very low density. 
Under these assumptions,  the population MC
 is null recurrent.
 Thus, the expected time of rarity $E_{X_0} T_{\e}$ is infinite, yet
any initially rare population will escape from rarity
 in finite time with probability one: $P_{X_0}(T_{\e}<\infty)=1$.
     We give a rigorous
proof that the law of the time spent in rarity, i.e.
$P_{X_0}(T_\epsilon = n)$, is
approximated by a power law with exponent $-3/2$ as $n\to \infty$.
 In contrast to the case $\l > 1$, this law is universal, in the sense that
it does not depend upon the distribution of the environmental 
noise. This behavior is the same as in the absence of 
population regulation (i.e. $a = 0$), when the population 
growth is driven by environmental noise only.

\begin{theo}
\label{th4}
  
   Let $f(x)=1+o(|\ln x|^{-1-\delta_2})$ as $x\to 0$ 
     and assume $E[|Y|^{2+\delta_1}]<\infty$
for some $\delta_1,\delta_2>0$.

      Then the MC is null recurrent.  This means that
  for any compact subset $A=[a,b]\subset ]0,\infty[$ and any $X_0\not\in A$, 
  $P_{X_0}(T_{A^c}<\infty)=1$ but $E_{X_0} T_{A^c}=\infty$.

  For any $p<1/2$ and for any $\delta>0$
  there exists $\epsilon(\delta)>0$ such that for all
  $X_0 <\e$
 \begin{equation}
 \label{b1}
 E_{X_0} T_{\e}^{p}=O(|\ln X_0|^{2p+ \delta}), \ \ X_0 \to 0.
 \end{equation}
   For any $p>1/2$, any $\epsilon>0$  and any $X_0<\epsilon$
\begin{equation}
\label{b2}
 E_{X_0} T_{\e}^{p}=\infty.
\end{equation}
  In particular for any $\epsilon>0$ small enough, any $X_0<\epsilon$
   and  any $\eta>0$ there exists a constant $C(X_0, \epsilon, \eta)>0$
 such that
  for all $n\geq 1$
\begin{equation}
\label{b3}
 C(X_0,\e,  \eta) n^{-3/2-\eta}\leq P_{X_0}(T_\e=n)
 \leq  C(X_0,\e, \eta) n^{-3/2+\eta}.
 \end{equation}
\end{theo}

 Compared to \cite{FH}
  we show that  null recurrence holds  also
  under assumption $Ee^{\a Y}=\infty$ for all $\a>0$.
   All estimates given by 
 the theorem  are new. In particular
   we recover  the predicted exponent $-3/2$ for $P_{X_0}(T_\e=n)$, 
or equivalently  $1/2$ for $P_{X_0}(T_\e>n)$ (cf. \cite{FC}).
Thus, in such `neutral'
populations, the time spent in  rarity  is much
longer than in commonness.
The tail of this time is universal,
i.e. it does not depend on the tail of the environmental noise
$Y$ as soon as 
$Y$
 has  more than a second moment, and it is the same 
as for  $f(x)\equiv 1$ provided 
the growth rate $f(x)$ departs from one sufficiently slowly as density
increases from zero.

Ferriere and
Cazelles \cite{FC} observed in their simulations
that the exponent $-3/2$ arises for a suitable limit
where $\l$ goes to one and $\e$ approaches zero.  We here recover this
result in the particular case where $\l$ is exactly one.
 The inequality (\ref{eqc}) also shows that
if $\sup_{x\le \e} \ln f(x)$ goes
to zero with $\e$, the exponent has to be greater than
$1$. It is  not yet clear how the $-3/2$
exponent could be recovered by taking such a limit.

\section{Two-species interaction}

The question of how these results extend to communities of
 interacting species opens a promising  avenue for future research. As
a first step in this direction,
 we report that the stochastic dynamics of two 
competing species with positive intrinsic growth can be
transient -- a result contrasting with the monospecific
case in which recurrence is guaranteed. 
Here, transience
is understood
in the sense that the two-species
community will not visit all  possible 
ranges of densities  with positive probability.
In fact, for some choices of the parameters,
we show that there exists a parameter
$\d>0$ such that
$X^1_n(X^2_n)^{-\d}$ is bounded below for all times 
 with positive probability provided this inequality
holds at time $n=0$.  However, the population
$X^1_n$ will undergo phases of commonness and
rarity in finite time. Although the population
density $X^2$ is somehow `dominated' by $X^1$,
we could not resolve whether the population density $X^2$ 
was limited by some upper-boundary 
at all times with positive probability or not.

 A formal result can be stated by considering the
 two-dimensional Ricker model with parameters $a_{11},
  a_{12}, a_{21}, a_{22}>0$
\begin{eqnarray}
\label{twR}
 X_{n+1}^1&=&X_n^1 e^{-a_{11}X_n^1-a_{12}X_n^2+r_1+ Y_n^1}\\
 X_{n+1}^2&=&X_n^2 e^{-a_{21}X_n^1-a_{22}X_n^2+r_2+ Y_n^2}\nonumber \\
\end{eqnarray}
   on $]0,\infty[^{2}$ where $( Y_n^1)_{n\in\N},(Y_n^2)_{n\in\N}$ are two independent  sequences of
 i.i.d.\ random variables.

\begin{theo}
\label{hsh}
   Assume that $E[(Y^1)^2 +(Y^2)^2]<\infty$,
        $r_1, r_2>0$,
    $r_1\ne r_2$ and
    one of the following  conditions is fullfiled:
\begin{enumerate}
\item $r_1 a_{21}-r_2 a_{11}>0$, $r_1 a_{22}-r_2 a_{12}>0$

\item  $r_1 a_{21}-r_2 a_{11}>0$, $r_1 a_{22}-r_2 a_{12}<0$

\item  $r_1 a_{21}-r_2 a_{11}<0$, $r_1 a_{22}-r_2 a_{12}<0$

\end{enumerate}
Then the MC is transient.

  Furthermore, we can describe the following ways of escape
   to infinity.

1. Assume that either condition (1) holds or
   $r_1 a_{21}-r_2 a_{11}>0$  together
  with the assumption $r_2< r_1$.
 Let us fix $\e>0$ such that $(r_1-\epsilon) a_{21}-r_2 a_{11}>0$.
   Let $M>0$. Let
\begin{equation}
\label{tMM}
\tau^{M}=\inf\{n: r_2 \ln X_n^1-(r_1-\epsilon)\ln X_n^2<M\}
\end{equation}
 if $r_2 \ln X_0^1-(r_1-\epsilon) \ln X_0^2>M$. 
Then
\begin{equation}
\label{fp1}
P(\tau^M=\infty)>0.
\end{equation}
Let $T_A(X^1)=\inf\{n \ge 0 ; X^1_n\in A^c\}$, $i=1,2$.
Then, for any $L,\e >0$,  
$$ E_{X_0}[T_L(X^1)]<\infty,\ E_{X_0}[T_\e(X^1)\mid \tau^M=\infty]<\infty,$$
\begin{equation}
\label{qmq1}
P_{X_0}(T_{[\e,L]^c}(X^1)\wedge \tau^M<\infty)=1.\end{equation}

2. Assume that either condition (3) holds or
   $r_1 a_{22}-r_2 a_{12}<0$  together
  with the assumption $r_2\geq r_1$.
 Let us fix $\e>0$ such that $(r_1+\epsilon) a_{22}- r_2 a_{12}<0$.
   Let $M>0$. Let
\begin{equation}
\label{tMM1}
\tau^{M}=\inf\{n:  r_2 \ln X_n^1-(r_1+\epsilon) \ln X_n^2>M\}
\end{equation}
 if $r_2 \ln X_0^1-(r_1+\epsilon) \ln X_0^2< M$. Then
\begin{equation}
\label{fp2}
P(\tau^M=\infty)>0.
\end{equation}
Let $T_A(X^2)=\inf\{n \ge 0 ; X^2_n\in A^c\}$.
Then, for any $L,\e\ge 0$,
$$E_{X_0}[T_L(X^2)]<\infty,\ E_{X_0}[T_\e(X^2)\mid \tau^M=\infty ]<\infty,$$
\begin{equation} 
\label{qmq2}
P_{X_0}(T_{[\e,L]^c}(X^2)\wedge \tau^M <\infty)=1.
\end{equation} 
\end{theo}
We see that the  two-population
MC can be transient even for some
  $r_1>0, r_2>0$ while in the 
one-dimensional case it is always recurrent
  for $r>0$.
  We conjecture that in the remaining case $r_1a_{21}-r_2 a_{11}<0$
  and $r_1 a_{22}-r_2 a_{12}>0$ the chain is recurrent,
as it can be suspected 
  from the vector field  of mean drifts, or from the isocline
structure of the 
  associated deterministic dynamical system.

\section{Discussion}

The importance of identifying the biologically 
relevant timescales of
ecological models, rather than focusing on their asymptotics,
 albeit
emphasized only recently [Hastings \cite{H}], has long been 
recognized in
some areas of population biology. In the study of epidemic dynamics, 
for
example, the relevant timescale can be the scale of a single 
outbreak, or
longer if the question of interest is the timing
 between epidemics
[Finkenstadt and Grenfell \cite{FG}]. Another example
 arises in the study
of plankton community dynamics in temperate regions, 
where  seasonality can
reduce the relevant timescale to less than one year, 
making the
asymptotic behavior of models meaningless [Huisman and 
Weissing, 
\cite{HW1}, \cite{HW2}]. 

Environmental stochasticity 
can play a crucial
role in determining the relevant timescales of a population's dynamics.
For example, adding stochasticity to a metapopulation
 models could turn
convergence to a stable focus equilibrium into a kind 
of dynamics that
would likely be identified as a noisy cycle by empiricists
 [Gurney et al.
\cite{NGP}], implying that the deterministic 
equilibrium provides limited
biological insight; instead, the  period of the apparent cycle
unravels a relevant
timescale to interpret fluctuations in such a system.
 Even small noise can
induce rapid fluctuations of large amplitude 
in population models that
would otherwise converge to stable equilibria 
[Higgins et al.
\cite{Hall}], e.g. by stabilizing population trajectories
 around chaotic
repellors [Rand and Wilson \cite{RW}].

Gyllenberg et al. \cite{GS} were the 
first to establish conditions under
which recurrence and null recurrence occur in the 
Ricker model studied in this
paper; they also addressed the case where the noise affects
 our parameter
$a$, which represents a random perturbation of 
the carrying capacity
of the habitat. Fagerholm and Hognas \cite{FH} extended
 their results by
considering the case where our parameter $a$ is also affected
 by i.i.d.
random perturbations. Our analysis offers a natural 
continuation of these
previous studies that (i) assumes less stringent conditions 
on the moments
of the random variables, (ii) provides precise and rigorous 
estimates of
the laws of several characteristic 
times, (iii) 
paves the way for tackling similar issues in 
 multispecific communities. 

Extension (i) is
interesting as it broadens the scope of the model to the
 case where noise
reflects random heterogeneity. Extension (ii) has several important
implications. Building on previous heuristics and numerical
studies, it yields a rigorous basis for the phenomenon of on-off
intermittency. The notion of on-off intermittency refers to a 
certain class of
burst-and-crash dynamics;  it was
first described in physical systems 
(Heagy et al. \cite{HPH}) and
later found to be relevant in the ecological context 
(Ferriere and Cazelles
\cite{FC}).
We have shown (cf. Theorem \ref{th02})
that for unbounded noise
the distribution of
time  spent in a medium abundance 
has exponentially decaying tail. 
Therefore, the population spends 
most of the time fluctuating 
in states of  commonness or rarity. The time 
spent in these states can be influenced by 
 the tail
of the environmental  noise when $\l\ge 1$.
If $\l>1$, the distribution of the exit time of these 
extreme states can have heavy tail (c.f Theorem \ref{th12})
even though the distributions of the time spent in
either
commonness or rarity  have exponential time.
If $\l=1$, the population spends most
of the time in rarity with
occasional short outbursts (cf. Theorem \ref{th4}).
For such neutral species, the distribution of rarity
 times can be approximated by
a power law with exponent $-3/2$; this approximation 
is universal in the
sense that it does not depend upon the law of the noise,
provided the latter  has finite second moment.
 Ferriere and
Cazelles \cite{FC} found a good fit of this 
approximation to real data on
fish population dynamics.

Many previous mathematical studies of population 
models based on
stochastic difference equations trace back to the
 seminal work of Ellner
\cite{E}. The first detailed study of a noisy version 
of the Ricker model
similar to ours was carried out by Schaffer et al. \cite{SEK}.
 Kornadt et
al. \cite{KLL} have considered a different model of environmental
stochasticity and investigated how stability of the non-trivial
equilibrium and  period-two cycles of the deterministic Ricker map
 was
affected. Similar questions were raised by Sun and Yang \cite{SY}
 for a
model of noise reflecting random immigration. Hognas \cite{H97} 
studied
the quasi-stationary distribution of the MC on a countable state 
space
(branching process) that describes a population regulated by the 
Ricker
mechanism with $\l > 1$. In Athreya and Dai \cite{at}, the state space
 is
continuous but the regulatory mechanism is given by the simple 
logistic
map; they study convergence and the existence of a stationary 
measure.
Ramanan and Zeitouni \cite{RZ} considered small noise operating
 additively
on an iterative map of a compact interval, like the logistic map. 
In their
case, extinction (i.e. transience in this study) was possible
 even with $\l>1$, and a tension arose between this effect of stochasticity 
and the
deterministic dynamics; for small noise, this tension resulted 
into large
extinction times.
There have been a number of large deviation analysis conducted 
on
discrete-time models  relevant for ecological applications.
 Morrow and Sawyer \cite{MS}   studied 
the large deviation tail of the exit time 
from a neighborhood of a stable fixed point 
of a MC in the limit where the noise vanishes 
as time goes to infinity; Kifer
 \cite{Ki1,Ki2} studied in a broad generality
the large deviation estimates for small perturbations
of dynamical systems.

Whether 
the universal power law
found for a single neutral population 
 exists for  multi-specific communities
  warrants
further investigation. From an ecological viewpoint, it would be
interesting to elucidate the potential relationship between 
the exponent of
such putative laws, and the structure of the network 
of interactions---how
many species involved, who interact with whom, and how? 
The goal then might
be to find some signature of a community's structure in
 the time series of
those neutral or quasi neutral
 species (i.e. having a $\l$ close to one)
 exhibiting 
on-off dynamics. Another 
desirable extension should aim at incorporating 
 temporal autocorrelation in the noise.
Temporal autocorrelation in random environments 
can have dramatic consequences on the dynamics of 
density-dependent  populations
(Kaitala et al. \cite{KYRL},
Ripa and Lundberg \cite{RL}).
Yet, the effect of autocorrelation at the characteristic timescales of
rarity and commonness remains completely unknown.

\section{Appendix :  proofs of the theorems}
The proofs of all our statements rely on
 martingale techniques.
 
\medskip

\noindent{\it Notation.} Everywhere below
   we denote by $E_{X_0, \ldots, X_i}$ the conditional expectation
 over the sigma-field ${\cal F}_i= {\cal F}(X_0,\ldots, X_i)$.

\subsection{Proof of  Theorem~\ref{th0}. }

    First, we need to get
 an a priori estimate for $E_{X_0} T_M$ (\ref{f1n}).
For any 
  $X_0$  large enough, define
$$\ln \bar X_n=\ln X_0 + \sum_{i=1}^{n}
\Big((\ln X_i-\ln X_{i-1})\vee (-1)\Big).$$
 Due to the assumption (\ref{ara})  on the function $f(x)$, 
   we may assume below that 
 $f(x)\le e^{-ax+r+\e}$ for all $x\ge M$
and some $\e>0$.
  Then
\begin{eqnarray*}
\lefteqn{E[ \ln \bar X_n- \ln \bar X_{n-1} \mid  X_{n-1}> x]}\\
&\le&E[(-aX_{n-1}+r+Y_n+\e)\vee (-1)\mid  X_{n-1}> x]\\
&\le &
(-ax+r+\e+1)P(Y>-1+a  X_{n-1} -r\mid X_{n-1}>x)
-1\\
&&{}+ E[Y1_{\{Y>-1+a X_{n-1}-r-\e\}} \mid X_{n-1}>x ]\\
&\le& -1+ E[Y1_{\{Y>-1+a x-r-\e\}}]=-1+o(1),\ \ x\to \infty
\end{eqnarray*} 
where we used  the fact
$E[Y1_{Y>-1+a x-r-\e}]\to 0$  as $x\to \infty$  due to the
 assumption
 $E|Y|<\infty$.
Thus for all $x$ large enough
  $$ E[\ln \bar X_n-\ln \bar X_{n-1} \mid X_{n-1}> x] \leq -1/2 $$
  for all $n\geq 1$.
It follows that
  $$E[\ln \bar X _{n \we T_M}
   \mid {\cal F}_{n-1}] \leq \ln \bar X_{n-1 \we T_M}
     -\frac{1}{2} 1_{\{T_M>n-1\} },$$
 hence
   $$E_{X_0}[\ln \bar X_{n \we T_M}]  \leq E_{X_0} [\ln \bar X_{n-1
 \we T_M}]
     -\frac{1}{2} P_{X_0}(T_M >n-1),$$
    and therefore
  $$E_{X_0}[\ln \bar X_{n \we T_M}]  \leq
     -\frac{1}{2} \sum_{i=1}^{n} P_{X_0}(T_M>i-1)+\ln X_0.$$
   But $\ln \bar X_{n} \geq \ln X_n$
     and  $\ln \bar X_n\geq \ln \bar X_{n-1}-1\geq \ln X_{n-1}-1$
so that $\ln \bar X_{n\we T_M}\ge \ln M-1$
  by the definition of $T_M$ and $\ln \bar X_n$. Then
$$ \ln M -1 \leq -\frac{1}{2}  \sum_{i=1}^{n} P_{X_0}(T_M>i-1)+\ln X_0 $$
   for all $n\geq 1$, whence
\begin{equation}
\label{f1n}
    E_{X_0} T_M= \sum_{i=1}^{\infty} P_{X_0}(T_M>i-1)\leq 2(\ln
    X_0-\ln M+1).
\end{equation}
Note that this inequality was obtained under the mere 
assumption that $f(x)\le e^{-ax+r+\e}$ for all $x\geq M$ and some $\e>0$. 
The next result requires the more precise  estimate
that for $x\ge M$  with $M$ fixed large enough
$$ e^{-ax+r-\e}\le f(x)\le e^{-ax +r+\e}$$
for some $\e\in (0,{1\over 2}r)$.
We have the following renewal equation:
\begin{equation}
\label{ws0}
E_{X_0}T_M = P_{X_0}(T_M=1)+E_{X_0}[ 1_{X_1>M} E_{X_1}[T_M+1]].
\end{equation}
The first term in this sum is the probability that
 the set $]0, M]$ is reached
 in one step.  Since  $f(x)\le e^{-ax+r+\e}$ for  $x\ge M$,
we find that
\begin{eqnarray}
\label{zz}
P_{X_0}(T_M=1)&=&P_{X_0}(X_1\le M)\ge P_{X_0}(X_0e^{-aX_0+r+\e+Y}
\le M)\nonumber\\
&= &P(Y<\ln M-\ln X_0 +aX_0 -r-\e).
\end{eqnarray}
 By the assumption  $E|Y|^{1+\delta}<\infty$
 we have  $P(Y>t) \leq C t^{-1-\delta}$ with some $C>0$ and all $t>0$, 
which yields  when $\ln M-\ln X_0 +aX_0 -r-\e\ge 0$
\begin{equation}\label{ws00}
P_{X_0}(T_M=1)\ge 1-C (\ln M-\ln X_0 +aX_0 -r-\e)^{-1-\delta}.
\end{equation}
  Thus 
\begin{equation}
\label{znn} 
P_{X_0} (T_M=1) \to 1, \ \ X_0 \to \infty.
\end{equation}
The remaining term is bounded according to  (\ref{f1n})
by
 $$E_{X_0}[ 1_{X_1>M} E_{X_1}[T_M+1]]\le E_{X_0}[
 1_{X_1>M } (1+ 2(\ln
    X_1-\ln M+1)) ]$$
$$\qquad\qquad\le E_Y( 1_{Y>\ln M-\ln X_0 +aX_0 -r-\e}
(1+2(Y+\ln X_0 -aX_0 +r+\e)
    -2\ln M+2)).$$ 
Again the assumption $E|Y|^{1+\delta}<\infty$
and Chebychev's inequality
imply, with $A=\ln M-\ln X_0 +aX_0 -r-\e$,  that
\begin{eqnarray}
E[ 1_{Y>\ln M-\ln X_0 +aX_0 -r-\e}
\lefteqn{(1+2(Y+\ln X_0 -aX_0 +r+\e)
    -2\ln M+2)]}\nonumber\\
&=& E[1_{Y\ge A}(3+2(Y-A))]
\le A^{-\d}E[|Y|^\d(3+2|Y|)]\nonumber \\
&\le & C'(\ln M- \ln X_0 +a X_0 -r-\e)^{-\delta}
\label{ws000}
\nonumber
\end{eqnarray}
 with some constant $ C'>0$
   for all $X_0$ large enough.
Hence 
\begin{equation}
\label{zlz} 
E_{X_0}[ 1_{X_1>M} E_{X_1}[T_M+1]]\to 0,  \ \ X_0\to \infty.
\end{equation} 
Combining (\ref{ws0}), (\ref{znn})
and (\ref{zlz}) we deduce that
  $E_{X_0}T_M \to 1$ as $X_0 \to \infty$.
  Furthermore by (\ref{ws00}) 
$$
P_{X_0} (T_M >1) \leq C_1 X_0^{-1-\delta}
$$
 with some constant  $C_1>0$.  Then $T_M \to 1$ a.s. as $X_0 \to \infty$
 by Borel-Cantelli lemma.

 Let us now prove the estimate (\ref{e1}). By Chebyshev's inequality
  for all $\a>0$
\begin{eqnarray}
\lefteqn{P_{X_0}(T_M>n)=P_{X_0}(\ln X_i>\ln M \ \
   \forall i=1,\ldots, n)}\nonumber\\
&=& P_{X_0} (e^{\a (\ln X_n-\ln X_0)} \geq e^{\a (\ln M-\ln  X_0)}, \
     X_i> M \ \ \forall i=1, \ldots, n-1)\nonumber\\
&\leq & e^{-\a(\ln M-\ln X_0)} E [e^{\a (\ln X_n-\ln X_0)}
    1_{\{  X_i> M \forall i=1, \ldots, n-1\}}]\nonumber\\
&=& (X_0/M)^{\a}
      E_{X_0}\Big[ e^{\a(\ln X_1-\ln X_0)}1_{\{ X_1>M\}}
        E_{X_0, X_1}\Big[ e^{\a(\ln X_2-\ln X_1)}
           1_{\{X_2>M\}} \nonumber\\
&& \cdots 1_{\{X_{n-2}>M\}}
       E_{X_0, \ldots, X_{n-2}} \Big[e^{\a(\ln X_{n-1}-\ln X_{n-2})}
      1_{\{X_{n-1}>M\}}\nonumber\\
&&\ \ \ \ \ {}          E_{X_0, \ldots, X_{n-1}}
 \Big[e^{\a(\ln X_{n}-\ln X_{n-1})}
    \Big]\Big]\cdots \Big]\Big].\label{e11}
\end{eqnarray}
       Note that for any $i=0, 1, \ldots, n-1$
 $$ 1_{\{X_i>M\}}
     E_{X_0, \ldots, X_{i}} e^{\a (\ln X_{i+1}-\ln X_{i})}
   \leq \sup_{x\ge M} f(x)^\a  E e^{\a Y}.$$ 
   Applying this inequality subsequently
  for $i=n-1, n-2, \ldots, 0$  to the right-hand
  side of (\ref{e1}) we  bound its left-hand side
  by $(X_0/M)^\a \kappa(\a)^{n}$ where
     $\kappa(\a)=\sup_{x\ge M} f(x)^\a E e^{\a Y}$.

   Our final task is to optimize this bound over $\a>0$.
  The function $\kappa(\a)$ exists  at least for $\a\in [0, \a_0]$ as
  $E e^{\a_0 Y}<\infty$.
     Furthermore $\kappa(0)=1$ 
   and $\kappa'(0)=\sup_{x\geq M} \ln f(x)<0$, from where 
  $\kappa(\a)<1$ for $\a>0$ small enough.
  Thus  $\kappa=\inf_{\a>0} \kappa(\a)<1$ and
  (\ref{e1}) follows.

\subsection{Proof of Theorem \ref{th02}}
The proof of this theorem is straightforward
since
\begin{eqnarray*}
P_{X_0}(T_{[\e,M]}\ge n)
&=& P_{X_0}\left( \cap_{1\le k\le n-2}
\{X_k f(X_k) e^{Y_k}\in [\e,M], X_k\in [\e,M]\}\right)\\
&\le&P_{X_0}\left( \cap_{1\le k\le n-2}
\{ K_{\rm min}(\e,M)
+\e\le Y_k \le K_{\rm max}(\e,M)
+M\}\right)\\
&=& \kappa(\e,M)^{n-2}\\
\end{eqnarray*}
with $K_{\rm min}(\e,M)=-\log\max_{x\in [\e,M]} (x f(x))$,
$ K_{\rm max}(\e,M)=
-\log\min_{x\in [\e,M]} (x f(x))$, and 
$$\kappa(\e,M)=P_{X_0}\left( K_{\rm min}(\e,M)
+\e\le Y \le K_{\rm max}(\e,M)
+M
\right)<1$$
since we assumed that $Y$ is
not compactly supported
and $0<\min_{x\in [\e,M]} (x f(x))\le \max_{x\in [\e,M]} (x f(x))
<\infty$.

\subsection{Proof of  Theorems \ref{th11}
and \ref{th12}.}

 Let us start by the proofs of (\ref{f2}) and (\ref{f3}).
    Assume that $\e$ is  fixed small
   enough. We fix $\e_0>\e$ and  introduce $\tilde X_n=X_n\wedge \e_0$.
  Then
\begin{eqnarray} 
\lefteqn{ E(\ln \tilde X_n   -\ln \tilde X_{n-1}\mid   X_{n-1}= x) }\nonumber\\
&=& 
  E[(-\ln f(x)+Y)1_{\{Y<-\ln (x f(x)) +\ln \e_0\}}]\label{skk}\\
&&{}+ E[(\ln \e_0-\ln x)1_{\{Y>-\ln (xf(x)) +\ln \e_0\}}]=\ln f(0) +o(1)\geq
\ln f(0)/2,\ \ x\to 0.\nonumber
\end{eqnarray}
  Here we used the facts that $\ln f(0)>0$ and, due to $EY=0$ and
   $E|Y|^{1+\delta}<\infty$,
$$E[Y 1_{\{Y<-\ln (xf(x)) +\ln \e_0\}}]=
  E[Y1_{\{Y\geq -\ln (xf(x)) +\ln \e_0\}}]\to 0,\ \ x\to 0,$$
 and by Chebyshev's inequality
 $$P(Y\geq -\ln (xf(x))+\ln \e_0) \leq \frac{ E |Y|^{1+\delta}}
   {(-\ln (xf(x)) +\ln \e_0)^{1+\delta}}=o((-\ln x)^{-1}),$$
 as  $x\to 0.$ 
 Thus $\ln \tilde X_{n\we T_\e }$ is a negative submartingale and
for $\e$ sufficiently small
$$E_{X_0}[\ln \tilde
  X_{n\we T_\e}]\ge E_{X_0}[\ln \tilde X_{n-1\we T_\e}]
+\frac{\ln f(0)}{2}P_{X_0}(T_\e \ge n-1)$$
resulting, for all $n\ge 0$, with
\begin{equation}
\label{zzf}
\ln \e_0 \geq E_{X_0}[\ln \tilde X_{n\we T_{\e}}]
  \geq  \ln X_0+\frac{\ln f(0)}{2}\sum_{p=0}^{n-1}P_{X_0}(T_\e \ge p).
\end{equation}
Hence, $T_\e<\infty$ almost
surely as a consequence of  Borel Cantelli's lemma. Moreover
\begin{equation}\label{gozo}
\frac{\ln f(0)}{2}E_{X_0}[T_{\e}]\le \ln \e_0-\ln X_0
\end{equation}
  which proves (\ref{f2}).

To get the lower bound on $E_{X_0}[T_\e]$ claimed in (\ref{f3}), first of
all, starting from the  rough estimate
\begin{equation}
\label{rou}
 E[\ln X_n- \ln X_{n-1} \mid  X_{n-1}<  \e ]
 \leq E_Y[(d +Y)] =d 
\end{equation}
   with  $d=\sup_{x\in [0,\epsilon]} \ln
  f(x)>0$  and
 proceeding in the same way as for (\ref{zzf})  we
  obtain the upper  bound
\begin{equation}
\label{ret}
  E_{X_0} \ln X_{n \we T_\e} \leq d \sum_{i=1}^{n}P_{X_0}
   (T_\e>i-1)+\ln X_0.
\end{equation}
  Since $T_\e<\infty$ a.s. by (\ref{f2}),
    then $\ln X_{n\we T_\e} \to  \ln X_{T_\e}$
 a.s.
 Note also that
  $$|\ln X_{n \we T_\e}| \leq Z=\sum_{k=1}^{\infty}
   |\ln X_{k \we T_\e}- \ln X_{k-1 \we T_\e}|+ \ln X_0.$$
  We have
$$E (| \ln X_{k \we T_\e}- \ln X_{k-1 \we T_\e}|
 \mid {\cal F}_{k-1})$$
$$ =   E (| \ln X_{k \we T_\e}- \ln X_{k-1 \we T_\e}|
1_{T_\e>k-1} \mid {\cal F}_{k-1}) \leq (d +E|Y|)1_{T_\e >k-1}.$$
  Then $E | \ln X_{k \we T_\e}- \ln X_{k-1 \we T_\e}|< C P(T_\e >k-1)$
  with some constant $C>0$.
  We have already shown that $E_{X_0} T_\e=
    \sum_{i=0}^{\infty}P_{X_0}(T_\e>i)<\infty$,
  then $EZ$ is finite.
  Hence,
  the dominated convergence theorem applies to the sequence $\ln X_{n \wedge T_{\e}}$:
$$E_{X_0}  \ln X_{n \we T_\e} \to E_{X_0}  \ln X_{T_\e}.$$
   But from the  definition of $T_\e$
  we have 
$E_{X_0} \ln X_{T_\e}\geq \ln \e$.   Thus taking the limit as $n\to \infty$
in (\ref{ret}) we  get
  the lower bound for $E_{X_0} T_\e$ valid for all $X_0$ small enough:
\begin{equation}\label{gozo2}
(\ln \e-\ln X_0)/d \leq E_{X_0} T_\e \end{equation}  showing
(\ref{f3}).

The proof of  (\ref{e2}) is analogous to that of  (\ref{e1}). Namely,
\begin{eqnarray}
\lefteqn{P_{X_0}(T_\e>n)=P_{X_0}(-\ln X_i>-\ln \e \ \
   \forall i=1,\ldots, n)}\nonumber\\
&=& P_{X_0} (e^{\a (-\ln X_n+\ln X_0)} \geq e^{\a (-\ln \e+\ln  X_0)},
     X_i< \e,\ \  \forall i=1, \ldots, n-1)\nonumber\\
&\leq & (\e/X_0)^{\a} E [e^{\a (-\ln X_n+\ln X_0)}
    1_{\{  X_i< \e \ \forall i=1, \ldots, n-1\}}]\label{e21}
\end{eqnarray}
   for all $\a>0$.
   We rewrite the right-hand side of (\ref{e21})
  as a sequence of  conditional expectations like in (\ref{e11}).
   Each of them is bounded from above  by:
$$ 1_{\{X_i<\e\}}
     E_{X_0, \ldots, X_{i}} e^{\a(-\ln X_{i+1}+\ln X_{i})}$$
$$   \leq 1_{\{X_i<\e\}} f(X_i)^{-\a}
     E e^{-hY} \leq  e^{-\a\inf_{[0,\e]} \ln f}E[ e^{-\a
     Y}]=\rho(\a ).$$
     Then the right-hand side of  (\ref{e21})
  obeys the upper bound $(\e/X_0)^{\a } \rho(\a)^n$
 with $\rho(\a)$ finite  for all  $\a \geq 0$
  small enough, $\rho(0)=1$ and $\rho'(0)=-\inf_{[0,\e]} \ln f<0$ for $\e$
small enough.
 Hence $\inf_{\a>0} \rho(\a)<1$ and (\ref{e2}) is proved.

     We are now ready to prove the recurrence of the MC.
   We denote by $A$ the compact
  $A=[\e, M]$, then $T_{[\e, M]^c}=\inf\{n\geq 1: X_n \in A\}$.
  We fix $\e>0$ small enough and $M$ large enough 
     such that 
  $$P(Y>\ln M-\sup_{x\in [0, \epsilon]} \ln f(x) -\ln \e)=\gamma<1.$$
      We prove that for any current density $X_0\in A^c$
   the time  to reach  $[\e, M]$  satisfies
\begin{equation}
\label{recc}
T_{[\e, M]^c}<\infty \ \ \hbox{a.s. }
\end{equation}
  Let
  $T_{1,M}=\inf\{n\geq 0 : X_n <M\}$,
  $T_{1,\e}=\inf\{n > T_{1,M}: X_n> \e \},\ldots,$
   $T_{k,M}=\inf\{n>T_{k-1,\e}: X_n< M\}$,
    $T_{k,\e}=\inf\{n>T_{k,M}: X_n>\e\}$.
  Then by (\ref{f1}) and (\ref{f2})
\begin{equation}
\label{zeze}
 P_{X_0}(\cap_{k=1}^{\infty} \{T_{k,M}<\infty, T_{k,\e}<\infty\})=1.
\end{equation} But
$$ P_{X_0} (T_{[\e, M]^c}> T_{k,\e}) \leq
    P_{X_0}(X_{T_{1,\e}}> M, X_{T_{2,\e}} >M,
     \ldots, X_{T_{k,\e}}>M )$$
$$= E_{X_0}[1_{\{X_{T_{1,\e}}> M, X_{T_{2,\e}} >M,
     \ldots, X_{T_{k-1,\e}}>M\} } P_{X_{T_{k-1,\e}}}(X_{T_{k,\e}}>M)]$$
$$\le  E_{X_0}(1_{\{X_{T_{1,\e}}> M, X_{T_{2,\e}} >M,
     \ldots, X_{T_{k-1,\e}}>M\}})P(Y > \ln M -\sup_{x\in
     [0,\epsilon]}\ln f(x) -\ln\e )$$
$$\le \cdots \le 
 [P(Y >\ln  M -\sup_{x \in [0,\epsilon]} \ln f(x)  -\ln\e )]^{k}=
     \gamma^{k}  \to 0, \qquad k \to \infty,$$
 where $0\leq \gamma <1$.
It follows that
$$ P_{X_0}\big (\cup_{k=1}^{\infty} (T_{[\e, M]^c}< T_{k, \e})\big)
   =\lim_{k\to \infty} P_{X_0}(T_{[\e, M]^c}<T_{k,\e})=\lim_{k\to \infty}
 1-\gamma^k=1,$$
 implying together with (\ref{zeze}) that 
$T_{[\e, M]^c}<\infty$ a.s.

  Let us show  that  $E[e^Y]=\infty$ implies
     $E_{X_0}[T_{[\e,M]^c}]=\infty$
 which  proves that the recurrence of the MC is null.
 First of all, we remark  that
$E_{X_0} \ln X_{T_M}=  E_{X_0} (\ln X_{1}1_{\{T_M=1\}})
+ E_{X_0} (\ln X_{T_M}1_{\{T_M\geq 2\}}).$
   Here, due to the inequality
 $e^{-ax+r-\delta}\le f(x)\le e^{-ax+r+\delta}$ for $x\ge M$
   and some $\delta>0$, we have:
\begin{eqnarray}
 E_{X_0} [\ln X_{1}1_{\{T_M=1\}}]
&\ge& E_{X_0}[(\ln X_0-a X_0+r +Y_1-\delta)
    1_{\{Y_1<\ln M
   -r +a X_0-\ln X_0-\delta\}}]\nonumber\\
&>&-\infty\label{cc} 
\end{eqnarray}
 but
\begin{eqnarray}
 E_{X_0} (\ln X_{T_M}1_{\{T_M\geq 2\}}) &\leq&
    E_{X_0} (\ln X_2 1_{T_M \geq  2})\nonumber\\
&=& E_{X_0}[(\ln X_0-a X_0+2r+2\delta +Y_1\nonumber\\
&&-a X_0
   e^{-a X_0+r -\delta +Y_1} +Y_2)1_{\{Y_1> \ln M
   -r +a X_0-\ln X_0-\delta\}}]\nonumber\\
&=&-\infty, \label{ccc}\\
\nonumber
\end{eqnarray}
   whence $E_{X_0}[ \ln X_{T_M}]=-\infty$.
As we have  $E_{X_0} [\ln X_{T_M}1_{\{X_{T_M}\geq \e\}}]>-\infty$, then 
 $ E_{X_0} [\ln X_{T_M}1_{\{X_{T_M}< \e\}}]=-\infty$, that is 
\begin{equation}
\label{kk}
  E_{X_0} [|\ln X_{T_M}|1_{\{X_{T_M}< \e\}}]=+\infty.
\end{equation}
     Let us also observe that by  (\ref{f3})
 there exists a constant $C>0$ such that  for all $X_0<\epsilon$
\begin{equation}
\label{kk1}
 E_{X_0}T_{\e} \geq C |\ln X_0|.
\end{equation}
   We may now conclude from (\ref{kk}) and (\ref{kk1}) that
\begin{equation}
\label{sss}
 E_{X_0}[ T_{[\e, M]^c}]
  \geq  E_{X_0} E_{X_{T_M}}[T_{\e} 1_{\{X_{T_M}< \e\} }]
 \geq
  E_{X_0}[C |\ln X_{T^{M}}| 1_{\{X_{T_M}< \e\}}
]=+\infty.
\end{equation}
   The positive recurrence of the MC under the assumption
  $Ee^Y<\infty$   has been proved in \cite{FH}.
   The proof in the genereal case of $f(x)$ with $\ln f(0)>0$ 
   is completely analogous: it uses the same test function 
   $g(x)=x 1_{\{x\geq x_0\}}+a^{-1}|\ln x|1_{\{x<x_0\}}$ 
  (where $x_0$ is such that $ax_0=|\ln x_0|$)
  for which the expansions (\ref{P0}), (\ref{P1}) with $\a=1$ are valid. 
  Therefore we do not repeat it.  

    Finally we turn to the refinements
  claimed  in (\ref{f6}) and (\ref{f5}).
  Let us first prove (\ref{f5}). Observe that for any finite
$X_0\neq 0$, $E[e^{\a Y}]=+\infty$ implies that
\begin{equation}\label{lp}
E[1_{\{X_2\le 1\}} (\ln {1\over X_2})^\a]=+\infty.
\end{equation}
Indeed,
$$\ln X_2\le (\ln X_0-(aX_0-r)+2\delta+Y_0 +r +Y_1 ) -
(aX_0 e^{-(aX_0-r)-\delta})e^{Y_0}:= Z_1-Z_2$$
with $E[|Z_1|^\a]<\infty$
since $Y_1,Y_0$ have finite first
moment, but $E[(Z_2)^\a]=+\infty$
by our assumption. (\ref{lp}) follows readily.
We now prove that there exists $a_\a>0$ and $b_\a$ finite
such that for any $X_0<\epsilon$ 
$$ E_{X_0}[T_\e^\a]\ge a_\a 1_{ \{X_0\le \e\}}
(\ln{\e \over X_0})^\a+b_\a$$
This is enough to get the desired estimate
by starting at time $n=2$ and
translating $T_\e$ by $2$ as in the proof of the null recurrence.
Set $S_n=\ln X_n$ so that 
$$S_{n\wedge T_\e}= S_0+\sum_{k=0}^{{n\wedge T_\e}}\ln f(X_k) 
+\sum_{k=0}^{{n\wedge T_\e }} Y_k.$$
We define $M_n=\sum_{k=0}^{{n}} Y_k$.
  Since $d=\sup_{x\in [0, \e]} \ln f(x)>0$, then   
 $\ln f(X_k)$'s are non-negative in this sum and bounded 
   from above by $d$. 
   Then  for any time $n\ge 0$
$$S_{n\wedge T_\e}\le S_0+d ({n\wedge T_\e}) + M_{n\wedge T_\e}.$$
In particular, since on $T_\e\le n$, $S_{n\wedge
T_\e}\ge\ln \e$ 
we get
$$ \ln \e\le
 ( S_0+d(T_\e\wedge n) +M_{T_\e \wedge n})$$
so that for all $n\in\N$
$$
 1_{\{T_\e\le n\}}\ts [(\ln \e-S_0)\vee 0 ]
\le 1_{\{T_\e\le n\}}\ts [(d(T_\e \wedge n) +M_{T_\e \wedge n})\vee 0].$$
Integrating the power $\a>0$ of this
inequality and using that for any $a,b\in\R$,
any $\a\ge 0$,
$|a+b|^\a \le 2^{\a} (|a|^\a+|b|^\a)$, 
we obtain
$$
E[[(\ln {\e \over X_0} )\vee 0 ]^\a 1_{\{T_\e\le n\}}]
\le E [|d(T_\e\wedge n) +M_{T_{\e} \wedge n })|^\a 1_{\{T_\e\le n\}}]$$
$$
\le 2^{\a} d^\a E [(T_\e\wedge n)^\a ]
+2^{\a} E [ |M_{T_\e\wedge n}|^\a]  .$$
By Burkholder-Davis-Gundy inequality (see e.g \cite{RY},
Theorem~4.1 p. 160), since $(M_n)_{n\ge 0}$
is a martingale with $<M>_n =E[Y^2] n $,
we know
that for any $\a>0$, there exists a finite
constant $C_\a$ such that
$$ E [ |M_{T_\e\wedge n}|^\a]\le C_\a E[(T_\e\wedge n
)^{\a\over 2}]
\le C_\a( 1+ E[(T_\e\wedge n)^{\a}])$$
so that
we get
$$[(\ln {\e\over X_0} )\vee 0 ]^\a P(T_\e\le n)-2^{\a} C_\a
\le  2^{\a} (d^\a+C_\a) E[(T_\e\wedge n)^{\a}].$$
We can now let $n$ going to
infinity and use monotone convergence theorem
with the fact that $T_\e$ is almost surely finite
 to  conclude
$$[(\ln {\e\over X_0} )\vee 0 ]^\a-2^{\a}C_\a
\le 2^{\a} (d^\a+C_\a) E[T_\e^{\a}]$$
which finishes the proof
of  (\ref{f5}) with (\ref{lp}).

 Finally it remains to prove (\ref{f6}).
 Let us  introduce a  positive  function
$g(x)=|a^{-1}\ln x|^\a$ for
$x<x_0$ and $g(x)=x^\a$ for $x\ge x_0$
with $x_0$ chosen so that $g$ is continuous.
Here, $a$ is such that 
$\lim_{x\ra \infty}e^{ax -r }f(x)=1,$
  hence $a^{-1}x^{-1}\ln f(x)=-1+O(x^{-1})$
  as $x\to \infty$.
Then by the assumptions $EY=0$, $E|Y|<\infty$
  and $Ee^{\alpha Y}<\infty$
  we have the following asymptotic expansions
\begin{eqnarray}
\lefteqn{E[g(X_{n+1})-g(X_n)|X_n=x]}\nonumber \\
&=& E[ (|a^{-1}\ln (x f(x)) + a^{-1}Y|^\a - x^{\a})
1_{\{Y< -\ln f(x) -\ln (x/x_0)\}}]\nonumber \\
&&{}
+E[ (x^\a
f(x)^{\a}e^{\a Y}
 - x^\a)1_{\{
Y\geq -\ln f(x) -\ln (x/x_0)\}
}]\nonumber \\
&=&  x^{\a}(-\alpha)  \frac{\ln x}{a x}(1+o(1))= 
  -\a a^{-1}x^{\a-1}\ln x(1+o(1)), \ \ x\to \infty\label{P0}
\end{eqnarray}
 and
\begin{eqnarray}
\lefteqn{E[g(X_{n+1})-g(X_n)|X_n=x]}\nonumber \\
&=& E[ (|a^{-1}\ln (xf(x))  + a^{-1}Y|^\a - |a^{-1}\ln x|^{\a})
1_{\{Y< -\ln f(x) -\ln (x/x_0)\}}]\nonumber \\
&&{}
+E[ (x^\a f(x)^\a e^{\a Y}
 -|a^{-1}\ln x|^\a)1_{\{
Y\geq -\ln f(x) -\ln (x/x_0)\}
}]\nonumber\\
&=&  |a^{-1}\ln  x|^{\a} E[(|1+(\ln f(x)+Y)(\ln x)^{-1}|^{\a}-1)
 1_{\{Y< -\ln f(x) -\ln (x/x_0)\}}]\nonumber\\
&&{}+o(|\ln x|^{\alpha-1})\nonumber \\
&=& |a^{-1}\ln x|^{\a}
 \alpha (\ln f(x)) (\ln x)^{-1}(1 +o(1))\nonumber \\
&=&-\a a^{-\a}(\ln f(0))|\ln x|^{\a-1}(1+o(1)),\ \ x\to 0. \label{P1}
\end{eqnarray}
  Let us now  start from some $X_0 \notin [\e, M]$
   with $\e$ chosen small enough and $M$ large enough.
   Then
  $g(X_{n \wedge T_{[\e, M]^c}})$ is a positive supermartingale.
 Note that $|a^{-1}\ln x|^{\a -1}= g(x)^{1-{1\over\a}}$ for $x<x_0$
and
   $x^{\a-1}=g(x)^{1-{1\over\a}}$ for $x>x_0$. Thus
 with some constant $\beta>0$
\begin{equation}
\label{ty}
E[g(X_{n+1})-g(X_n) \mid X_n=x]\leq -\beta (g(x))^{1-{1\over\a}}
\ \ \forall x\not\in [\e, M],
\end{equation}
whence
\begin{equation}\label{ty2}
E_{X_0}[g(X_{(n+1)\wedge T_{[\e, M]^c}}
)-g(X_{n\wedge T_{[\e, M]^c}})]\le -\beta E_{X_0}[ g(X_{n})^{1-{1\over\a}}
1_{\{T_{[\e, M]^c}\ge n\}}].
\end{equation}
Assume $\a<1$ so that $x\ra x^{1-\a^{-1}}$ is decreasing. Then for any sequence $p_n> 0$,
$$E_{X_0}[ g(X_{n})^{1-{1\over\a}}1_{\{T_{[\e, M]^c}\ge n\}}]\ge
p_n^{1-{1\over\a}}E_{X_0}[1_{\{g(X_{n})\le p_n\}}
1_{\{T_{[\e, M]^c}\ge n\}}]$$
$$\ge p_n^{1-{1\over\a}}\Big(P_{X_0}(T_{[\e, M]^c}\ge
n)-{E_{X_0}(g(X_n)1_{\{T_{[\e, M]^c}\ge n\}})\over
p_n}\Big)$$
so that finally we deduce
\begin{eqnarray}
\lefteqn{\beta p_n^{1-{1\over\a}}P_{X_0}(T_{[\e, M]^c}\ge n)}\label{ty3}\\
&\le& E_{X_0}[g(X_{n\wedge T_{[\e, M]^c}})-g(X_{(n+1)\wedge T_{[\e, M]^c}}
)] +\beta p_n^{-{1\over\a}}E_{X_0}(g(X_n)1_{\{T_{[\e, M]^c}\ge n\}}).
\nonumber
\end{eqnarray}
Observe that by (\ref{ty})
$E_{X_0}(g(X_n)1_{\{T_{[\e, M]^c}\ge n\}})\leq  g(X_0)$.
Then summing the inequalities (\ref{ty3}) over  $n=0,1,\ldots,m$
yields
$$\beta \sum_{n=0}^m
 p_n^{1-{1\over\a}}P_{X_0}(
T_{[\e, M]^c}\ge n)\le g(X_{0})-E_{X_0}g(X_{m+1\wedge T_{[\e, M]^c}}
) +\beta g(X_0)\sum_{n=0}^m  p_n^{-{1\over\a}}$$
$$\le g(X_{0}) +\beta g(X_0)\sum_{n=0}^m  p_n^{-{1\over\a}}
 \ \ \forall m=1,2,\ldots.$$
Letting $m\to \infty$ we conclude that
$$\sum_{n\ge 0}  p_n^{-{1\over\a}}<\infty\Rightarrow
\sum_{n\ge 0}
 p_n^{1-{1\over\a}}P_{X_0}(T_{[\e, M]^c}\ge n)<\infty.$$
Now it remains to  take $p_n^{1\over \a} = n(\ln n)^{1+\zeta}$
   with $\zeta>0$ small enough
to get the right asymptotics of the tail of $T_{[\e, M]^c}$
since for any non negative  random variable $T$ and  any $a>0$
$$E[T^a]<\infty\quad\Leftrightarrow\quad
\sum_{n\ge 0} n^{a-1}P(T\ge n)<\infty.$$
If $\a\ge 1$, we set $Y_n= g(X_n)^{1\over 2\a}$
so that (\ref{ty}) becomes
\begin{equation}
\label{tyd}
E[Y_{n+1}^{2\a} -Y_{n}^{2\a} \mid Y_n=y
]\leq -\beta Y_{n}^{2\a-2}
\ \ y\ge L,
\end{equation}
with the special choice of $\e,M$
such that $L=g(\e)^{1\over 2\a}=g(M)^{1\over 2\a}$
(which we can always do up to take $\e$ smaller or $M$
larger).
Hence, by Theorem 1 of \cite{AIM},
we deduce, since 
$\{Y_{n}\ge L\}\subset \{T_{[\e,M]^c}>n\}$,
that for $\a\ge 1$, there exists $c<\infty$ such that
$$E[T_{[\e,M]^c}^\alpha]\le c g(X_0).$$

\bigskip

\subsection{ Proof of Theorem \ref{th4}}
 The results (\ref{b1}) and (\ref{b2})
 follow from the estimates on passage-time moments
 for nonnegative stochastic processes of \cite{AIM}.
 These are generalizations
of Lamperti's results for countable Markov chains \cite{L}.

   Consider the Markov chain on $[0, \infty[$ 
       $$Z_{n+1}= (Z_n-\ln f(e^{-Z_n})-Y_n)
   1_{\{Z_n-\ln f(e^{-Z})-Y_n>0, \; Z_n>0 \}}$$  that 
   starts at $Z_0>0$ and has $0$ as an absorption state. 
   Let $S_A=\inf\{n \geq 0: Z_n \in A^c\}$ the time to escape 
   from the subset $A \in [0, \infty)$. 
     We will show that for any $p<1/2$ and any $\delta_0>0$, there exists
   $\epsilon(\delta)>0$ such that  for any $0<\epsilon\leq \epsilon(\delta)$
\begin{equation} 
\label{zsz} 
 E_{Z_0} S_{]-\ln \e, \infty[}^{p}=O(Z_0^{2p+\delta}),\ \ Z_0\to
 \infty 
\end{equation}      
  and also that for any $p>1/2$ and any $\epsilon>0$ 
\begin{equation} 
\label{zsz1} 
 E_{Z_0} S_{]-\ln \e, \infty[}^{p}=\infty. 
\end{equation}
  Then the change of variables 
     $Z_n=-\ln X_n$ implies immediately the statements (\ref{b1}), 
      (\ref{b2}) of the theorem 
   as we have $T_\e=S_{]-\ln \e, \infty[}$ through this change. 
   It follows from (\ref{b1}) that    
  $T_{\e}<\infty$ a.s. for any $X_0<\epsilon$.   
   This fact entails the recurrence of the MC proceeding exactly by the same
     arguments as in the case $f(0)>1$.  
  The recurrence is null  by (\ref{b2}). 
   The statement (\ref{b3}) is a direct consequence of 
   (\ref{b1}) and (\ref{b2}). 
 
  To see (\ref{zsz}) and (\ref{zsz1}), we denote by 
$$\mu_r(z)=\left\lbrace
\begin{array}{l}
E(Z_{n+1}-Z_n|Z_n=z)\mbox{ if } r=1,\cr
E(|Z_{n+1}-Z_n|^r|Z_n=z)\mbox{ if } r\neq 1.\cr
\end{array}
\right.$$  
     By the assumptions $EY=0$, 
  $E|Y|^{2+\delta_1}<\infty$  we have
   $P(Y>z)=O(z^{-2-\delta_1})$ 
   $EY1_{Y<z}=O(z^{-1-\delta_1})$, 
   $EY^2 1_{Y<z}=O(z^{-\delta_1})$ as $z\to \infty$.
   By the asumtpion $f(x)=1+o(|\ln x|^{-1-\delta_2})$ 
   as $x\to 0$ with some $\delta_2>0$, we have
   $f(e^{-z})= 1+o(z^{-1-\delta_2})$ as $z\to \infty$. Then 
$$
\mu_1(z)= E((-\ln f(e^{-z}) -Y)1_{\{Y< z-\ln f(e^{-z})\}} 
     -z1_{\{Y>z-\ln f(e^{-z})\}})$$
$$= O(z^{-1-\delta_2})+
     O(z^{-1-\delta_1}),\ \ \ z\to \infty
$$  
$$
\mu_2(z)= E((-\ln f(e^{-z}) -Y)^2 1_{\{Y< z-\ln f(e^{-z})\}} 
     +z^2 1_{\{Y>z-\ln f(e^{-z})\}})$$
$$= O(z^{-2-2\delta_2})+O(z^{-2-\delta_1-\delta_2})+
     E Y^2 + O(z^{-\delta_1}),\ \ \ z\to
     \infty,
$$
and for any $2<r<2+\delta_2$ 
$$\mu_r(z)\leq 2^{r}(o(1)+|Y|^r)+z^{r} O(z^{-2-\delta_1})
   =o(z^{r-2}),\ \ z\to \infty.
$$   
  Then for any $p\in \R$  
$$ 2z\mu_1(z)+(2p-1)\mu_2(z)=(2p-1)E(Y^2)+O(z^{-\delta}) $$ 
   with some $\delta>0$.
  Consequently,
  Propositions 1 and 2 p.957 in \cite{AIM} apply to the MC $Z_n$ and  
  prove (\ref{zsz})  and (\ref{zsz1}) concluding the proof 
 of the theorem.

\subsection{Proof of Theorem~\ref{hsh}.}

   Clearly, the transience in cases (1), (2) and (3) follows from
 (\ref{fp1}) and (\ref{fp2}).
We restrict ourselves  to (\ref{fp1}),
   the second case being symmetrical.

   Let us construct the function
\begin{equation}
\label{cmc}
g(X^1, X^2)=(r_2 \ln X_n^1-(r_1-\epsilon) \ln X_n^2)^{-1} \vee d^{-1}
\end{equation}
  with some fixed $d<M$. We use the notation
   $\vec X= (X^1, X^2)$ for shortness.
  We shall prove that there exists $\beta_0>0$ such that
   for all $\beta\geq \beta_0$
 \begin{equation}
\label{zv}
 E(g(\vec X_{n+1})-g(\vec X_n) \mid g(\vec X_n)=\beta^{-1})\leq 0.
\end{equation}
  Then for all $n\geq 1$
$$ E g(\vec X_{n\wedge \tau^M}) \leq E g(\vec X_0)=
   (r_2\ln X_0^1-(r_1-\epsilon) \ln X_0^2)^{-1}.$$
   Assume now the contrary of (\ref{fp1});  $P(\tau^M<\infty)=1$.
 Then $ \vec X_{n\wedge \tau^M} \to X_{\tau^M}$ a.s. as $n\to \infty$.
  Then by Fatou's lemma
$$ E g(\vec X_{\tau^M}) \leq \lim_{n\to \infty}
        E g(\vec X_{n\wedge \tau^M}) \leq E g(\vec X_0)=
   (r_2\ln X_0^1-(r_1-\epsilon) \ln X_0^2)^{-1}.$$
This contradicts the definition of $\tau^M$
  as $g(\vec X_{n\wedge \tau^M})> M^{-1}$.

   The proof of (\ref{fp1}) is reduced now to (\ref{zv}).
Let us note that  under condition
    $r_2 \ln X_{n}^1- (r_1-\epsilon) \ln X_{n}^2=
   \beta$ we have:
\begin{equation}
\label{hhh}
r_2 \ln X_{n+1}^1- (r_1-\epsilon) \ln X_{n+1}^2
   =\beta +h_1 X_n^1+h_2 X_n^2 + h + Y_n
\end{equation}
  where
$$
h_1 =-r_2 a_{11}+ (r_1-\e)a_{12}>0,
 \ h_2= -r_2 a_{12}+(r_1-\e)a_{22},
  \ h=\e r_2>0,$$
$$\ Y_n= r_2 Y_n^1-(r_1-\e) Y_n^2.$$
 Our assumptions on $r_i$ and $a_{i,j}$ for $i,j=1,2$
  in this case ensure the following property:
     for any $\delta>0$ there exists $\beta_0>0$ such that
  {\sl for
   all } $X^1, X^2>0$  with
   $r_2 \ln X^1- (r_1-\epsilon) \ln X^2=\beta
      \geq \beta_0 $
\begin{equation}
\label{fasss}
h_1 X^1+ h_2 X^2 >-\delta.
\end{equation}
        In fact, if condition (1) holds, then also $h_2>0$.
    Then trivially  this quantity is non-negative.
     Otherwise, if $h_1>0$ but $h_2<0$ and $r_2<r_1$, we can rewrite
$$
h_1 X^1+ h_2 X^2 = h_1 e^{\ln X^1}+ h_2
e^{-\beta/r_2+[r_2/(r_1-\epsilon)-1]\ln X^1 +\ln X^1}
$$
If $X_n^1\to \infty$, the first term
    with positive coefficient $h_1>0$
     dominates  this sum due to
  the fact that  $r_2/(r_1-\epsilon)-1<0$. Hence,
  we have the desired  estimate for all $X^1>C$
  and all $\beta>0$.
It remains to chose $\beta_0(\delta)>0$ large enough
    to ensure  (\ref{fasss}) for $0<X^1<C$.
  For the rest of the proof we fix
    $\delta=h/2$. Thus
 \begin{equation}
 h_1 X_n^1+h_2 X_n^2+h\geq h/2>0 \ \ \forall \b\geq \b_0.
 \end{equation}
The left-hand side of (\ref{zv}) equals $I_{\b}^1+I_{\b}^2$
   where
\begin{eqnarray}
I_\b^1&=&E[(d^{-1}-\b^{-1}){\bf 1}_{\{Y_n< d-
\beta -h_1 X_n^1-h_2 X_n^2- h
  \}}]\nonumber\\
 &\leq & d^{-1}
P(Y_n< d-
\beta -h_1 X_n^1-h_2 X_n^2- h)\nonumber\\
 &\leq&
  \frac{ \kappa }
{(\b+h_1 X_n^1+h_2 X_n^2+h-d)^{2}\ln(\b+h_1 X_n^1+h_2 X_n^2+h-d)}
  \label{ehe}
\end{eqnarray}
  for some $\kappa>0$ as $E|Y|^2<\infty$. Here we used the fact that $\beta$
is big enough so that $\beta-d+(h/2)\ge 0$. 
  We have
\begin{eqnarray}
I_\b^2&=&E[((\b - h_1 X_n^1-h_2 X_n^2-h -Y_n)^{-1}
    -\b^{-1}){\bf 1}_{\{Y_n>d-
\beta -h_1 X_n^1-h_2 X_n^2- h
  \}}]\nonumber\\ 
&=& E\Big[\frac{-h-Y_n}{(h_1 X_n^1+h_2 X_n^2+h+\b)(h_1 X_n^1+h_2 X_n^2+\b)}
{\bf 1}_{\{Y_n>d-
\beta -h_1 X_n^1-h_2 X_n^2- h
  \}}\Big]\nonumber\\
&& {}- \frac{1}{(\b+h_1 X_n^1+h_2 X_n^2)^{2}}\nonumber\\
&&\ \ \ \ {}\times \Big(E[\frac{(\b+h_1 X_n^1+h_2 X_n^2)(h_1 X_n^1
    +h_2 X_n^2)}{\b} {\bf 1}_{\{Y_n>d-
\beta -h_1 X_n^1-h_2 X_n^2- h
  \}}]\nonumber\\
&&\ \ \ \ \ \ \ \ \ {}-\frac{\b+h_1 X_n^1+h_2 X_n^2}{\b +h_1 X_n^1
    +h_2 X_n^2 +h}\nonumber\\
&&\ \ \ \ \ \ \ \ \ \ \ \ \ \ {}\times  E\Big[\frac{Y_n(c+Y_n)}{\b+h_1 X_n^1+h_2
    X_n^2+h+Y_n}
{\bf 1}_{\{Y_n>d-
\beta -h_1 X_n^1-h_2 X_n^2- h
  \}}\Big]\Big).\nonumber
\end{eqnarray}
 One can bound it from above by
\begin{eqnarray}
I_{\b}^2&\leq &
 \frac{1}{(h_1 X_n^1+h_2 X_n^2+h+\b)(h_1 X_n^1+h_2 X_n^2+\b)}\\
&&{}\times \Big(-h+ h P(Y_n<d-\b-h/2)+
+ E [Y_n {\bf 1}_{\{Y_n<d-
\beta-h/2\}}]\Big)\nonumber\\
&&\ {}+ \frac{-1}{(\b+h_1 X_n^1+h_2 X_n^2)^{2}}\label{zd}\\
&&\ \ \   \times 
\Big(-h/2  -\frac{\b+h_1 X_n^1+h_2 X_n^2}{\b +h_1 X_n^1
    +h_2 X_n^2 +h}\nonumber\\
&&\ \ \ \ \ \ \ {}\times  E\Big[\frac{Y_n(c+Y_n)}{\b+h_1 X_n^1+h_2
    X_n^2+h+Y_n}
{\bf 1}_{\{Y_n>d-
\beta -h_1 X_n^1-h_2 X_n^2- h
  \}}\Big]\Big)\nonumber
\end{eqnarray}
  where in the first estimate we used the fact that $EY_n=0$.
 Here $P(Y_n<d-\b-h/2) \to 0$
   and $E Y_n {\bf 1}_{\{Y_n<d-\b-h/2\}} \to 0$
as $\b\to \infty$ since $E|Y|^2<\infty$.

       Finally the sequence
  of the r.v. $Z_{m} = \frac{Y_n(c+Y_n)}{m+Y_n}
{\bf 1}_{\{Y_n>d-m\}}$
 converges to zero a.s. as $m \to \infty$
   and $|Z_{m}|\leq d^{-1} |Y_n(c+Y_n)|$ where
     $E|Y_n(c+Y_n)|<\infty$.
Then by dominated convergence theorem
   $E Z_m \to 0$ as $m\to \infty$.
   Since by (\ref{fasss})  $\b+h_1 X_n^1+h_2 X_n^2+h \geq \beta+h/2$
 $$E\Big[\frac{Y_n(c+Y_n)}{\b+h_1 X_n^1+h_2
    X_n^2+h+Y_n}
{\bf 1}_{\{Y_n>d-
\beta -h_1 X_n^1-h_2 X_n^2- h
  \}}\Big]\to 0,\ \ \b \to \infty$$
  uniformly for  all $X_n^1,X_n^2>0$  with
  $r_2 X_n^1- (r_1-\e)X_n^2=\b$.
Combining these facts,
  we see that the sum of the estimates
  of (\ref{zd}) and (\ref{ehe}) equals
  $(\b+h_1 X_n^1+h_2 X_n^2 +h)^{-2}(o(1)-h+o(1)+h/2+o(1))$
     as $\b\to \infty$ uniformly for all
   $X_n^1, X_n^2$ with
       $ r_2 X_n^1-(r_1-\e)X_n^2=\beta$.
   Then (\ref{zv}) is satisfied and the
first part of the theorem is proved.
To prove (\ref{qmq1}), we notice
that
\begin{eqnarray}
X_{n+1}^1&\le& X_{n}^1
e^{-a_{11}X_{n}^1 +r_1 +Y^1_{n}}\label{lk1}
\\
X_{n+1}^1&\ge& X_{n}^1
e^{-a_{11}X_{n}^1-d(M) (X_{n}^1)^{\d} +r_1 +Y^1_{n}}\mbox{ on } \tau^M\ge n
\label{lk2}\\
\nonumber
\end{eqnarray}
with $d(M)=a_{12}e^{-{r_2\over r_1-\e} M}$
and $\d={r_2\over r_1-\e}$.
We therefore can prove  that $T_L(X^1)$
satisfies the analogue of the bound 
 (\ref{f1n}) of Theorem~\ref{th0}
by (\ref{lk1}) (note that we only used in the proof
  of (\ref{f1n})
the upper bound), whereas we can use 
(\ref{lk2}) to get bounds (\ref{f2}) and (\ref{f3})  
on $T_\e(X^1)$ under condition $\tau^M =\infty$.
In fact, we can apply the same arguments
than in the proofs of Theorems~\ref{th0}, ~\ref{th11}
and~\ref{th12} by considering the  $(\bar X_{n},
\tilde X_{n\wedge\tau_M})_{n\ge 0}$
instead of $(\bar X_{n},
\tilde X_{n})_{n\ge 0}$ yielding $E_{X_0}[ T_\e(X^1)\wedge
\tau_M]<\infty$ and $E_{X_0}[ T_L(X^1)]<\infty$
(note here that in the neighborhood of the
origin, the correction $d(M) (X_{n}^1)^\d$ is small).
     Then as in one-dimensional case,  
 we can derive that 
  $P_{X_0}(T_{[\e, L]^c}(X^1)\wedge\tau_M
<\infty)=1$. 

\bigskip 

\nn
{\bf Acknowledgment}
We are very grateful to O. Zeitouni
for pointing out some useful references.

\bibliographystyle{plain}
\bibliography{ricker}

\end{document}